\documentclass[12pt]{article}

\usepackage{amsmath}
\usepackage{amsthm}
\usepackage{amsfonts}
\usepackage{yfonts}
\usepackage{enumitem}
\usepackage{graphicx}
\usepackage[autolanguage]{numprint}

\usepackage{graphicx}
\usepackage{capt-of}

\newtheorem{theorem}{Theorem}

\usepackage[english]{babel}

\usepackage{tikz}

\begin{document}

\title{Orthogonal curvilinear coordinate systems and torsion-free sheaves \\ over reducible spectral curves
\thanks{The work was supported by the RSCF (grant 19-11-00044-P).}}
\author{	A.E. Mironov \thanks{Novosibirsk State University and Sobolev Institute of Mathematics, e-mail: mironov@math.nsc.ru},
A. Senninger \thanks{Novosibirsk State University, e-mail:
 alexandre.senninger@gmail.com},
 I.A. Taimanov \thanks{Novosibirsk State University and Sobolev Institute of Mathematics, e-mail: taimanov@math.nsc.ru}}
\date{}
\maketitle

\abstract{The methods of finite-gap integration are used to construct orthogonal curvilinear coordinate systems in the Euclidean space corresponding to sheaves of rank one without torsion over reducible singular spectral curves.}

\section{Introduction}
Orthogonal curvilinear coordinate systems in the Euclidean spaces play an important role in geometry, differential equations, mathematical physics and other branches of mathematics. The construction and classification of orthogonal curvilinear coordinates is a classical problem of differential geometry. This problem was studied by Gauss, Darboux, Bianchi, etc. In general, the construction of coordinates is reduced to finding coordinate functions $x^j=x^j(u)$ from the overdetermined system of differential equations:
$$
 \partial_{u^i}x^1\partial_{u^j}x^1+\dots+\partial_{u^i}x^n\partial_{u^j}x^n=0,\quad i\ne j.
$$

In \cite{Z} Zakharov initiated the application of the methods of soliton theory to this problem, using the dressing method. Shortly after this work, Krichever found a procedure for constructing such coordinates by finite-gap integration methods \cite{K1}. In Krichever's work, explicit theta-functional formulas for coordinates are given, which, as usual, include parameters which satisfy transcendental relations. In this case, theta functions correspond to smooth spectral curves. The problem of efficiency, which was posed by Novikov for theta-functional solutions of soliton equations, has not yet been investigated in this case.

At the same time, if we consider the limiting construction, when a spectral curve
(a two-dimensional Riemann surface) degenerates into a reducible curve with components of small genus, we can obtain explicit formulas. This is obviously the case when all irreducible components are rational curves, i.e. closed Riemann surfaces of genus zero.

Such a construction was obtained in \cite{MT1}, where, in particular, it was indicated how to use it to obtain the polar and spherical coordinates. In \cite{MT2} this approach is applied to the construction of new examples of Frobenius manifolds.

Other examples of using singular spectral curves to construct exact solutions are given in \cite{T}.


In \cite{MT1}, there was considered the case when the spectral curve is constructed from its normalization
$\Gamma_\mathrm{nm}$, a disjoint union of nonsingular curves, by identifying pairs of points on them.
\par
The Baker--Akhiezer function $\psi$ is given as a function on
$\Gamma_\mathrm{nm}$ satisfying additional relations of the form
\begin{equation}
\label{e1}
\psi(a) = \psi(b),
\end{equation}
where the points $a$ and $b$ are glued together and by the projection $\Gamma_\mathrm{nm}\to\Gamma$ are mapped into a double point on $\Gamma$.

\par\leavevmode\par

In this paper, we show that the results of \cite{MT1} can be generalized to reducible singular spectral curves equipped with torsion-free sheaves of rank one. In this case, the Baker-Akhiezer function, which determines the coordinates, is not a function, as in the case of a smooth spectral curve, but a section of such a sheaf. In terms of Baker-Akhiezer sections, this looks like a replacement of \eqref{e1} by a more general relation
\begin{equation}
\label{e2}
\psi(a) = \lambda \psi(b),
\end{equation}
where the parameter $\lambda$ can be different from one.

In \S 2 we present the Krichever construction and our main results. 
In \S 3 we present the spectral data for Darboux--Egorov metrics in the case of reducible spectral curves endowed with torsion-free sheaves of rank one. In \S 4 we give some concluding remarks on interesting problems related to these constructions.

\section{Baker--Akhiezer function and orthogonal curvilinear coordinates}

\subsection{The Krichever construction \cite{K1}}

Let $\Gamma$ be a smooth closed Riemann surface of genus $g$.
We assume that on $\Gamma$ there are given divisors
$$
P = P_1+\dots+P_n, \quad Q= Q_1+\dots+Q_n,
$$
$$
R = R_1 +\dots+ R_l, \quad \gamma = \gamma_1 +\dots+\gamma_{g +l-1},
$$
and local parameters $k_j^{-1}$ near the points $P_j$.
We assume that all points from these divisors are distinct.

There is a unique {\it$n$-point Baker-Akhiezer function $\psi(u,P)$} having the following properties:

1) in the neighborhood of $P_j$, it has the expansion:
$$\psi = e^{k_j u^j} \left( h_{j} (u) + \dfrac{\xi_{j} (u)}{k_j} +... \right),$$

2) on $\Gamma\setminus\{\cup P_j\}$ the function $\psi(u,P)$ has simple poles in $\gamma$,

3) $\psi(R_j) = d_j$ (not all equal to $0$).
\footnote{In \cite{K1} $d_j=1$ with $j=1,\dots,l$, however, in the more general case indicated by us, all the reasoning goes unchanged.}


The function $\psi$ can be explicitly expressed in terms of the theta function of the Jacobi variety of the spectral curve $\Gamma$.

\vskip 4mm

{\bf Theorem A\cite{K1}.}
{\sl Let on $\Gamma$ there exist a holomorphic involution $\sigma : \Gamma\rightarrow\Gamma$ with fixed points $P_j$ and $Q_j$, and $\sigma(k_j) = - k_j$.
Let there exist a meromorphic $1$-form $\Omega$ on $\Gamma$ such that its divisor of zeros and poles have the form:
$$(\Omega) = P_1 +\dots +P_n + \gamma + \sigma(\gamma)-Q_1 -\dots - Q_n - R - \sigma(R),$$
and
$$
{\rm Res}_{Q_1} \Omega = \dots ={\rm Res}_{Q_n} \Omega.
$$
Then
$$
\partial_{u^i} x^1 \partial_{u^j} x^1 +\dots+ \partial_{u^i} x^n \partial_{u^j} x^n = 0, \quad i \neq j,
$$
where
$$
x^j(u) = \psi (u,Q_j).
$$}

\vskip4mm

\subsection{Orthogonal curvilinear coordinates corresponding to the Baker--Akhiezer functions on singular spectral curves \cite{MT1}}


In \cite{MT1}, Krichever's construction of orthogonal curvilinear coordinates was generalized to the case of reducible spectral curves. Namely, let $\Gamma$ be a reducible spectral curve consisting of irreducible components $\Gamma_1,\dots,\Gamma_m$.
Let on $\Gamma$ be given divisors
$$
P =P_1+\dots+P_n, \quad Q= Q_1+\dots+Q_n,
$$
$$
R =R_1 +\dots+R_l, \quad \gamma= \gamma_1 +\dots+\gamma_{g_a + l-1},
$$
and $k_j^{-1}$ near the points $P_j$.

 In difference with \S 2.1, instead of the genus $g$ of the spectral curve $\Gamma$, we consider the arithmetic genus $g_a$ of the singular curve. It differs from the geometric genus (for smooth surfaces both kinds coincide), but it is included in the Riemann-Roch theorem for singular surfaces.

We define the function $\psi$ by the same conditions 1)--3) as in the case of a non-singular spectral curve.

The Baker--Akhiezer function $\psi$ on the reduced curve $\Gamma$ is given by the functions
$\psi_j$ on each component $\Gamma_j$. The union of the normalizations of these components forms the normalization of the spectral curve:
$$
\Gamma_\mathrm{nm} = \cup \Gamma_{j,\mathrm{nm}}
$$
(the facts on singular spectral curves necessary for finite-gap integration can be found in \cite{Serre}, see also \cite{T}).

Let us consider one particular but important case. Let there be a double point on $\Gamma$,
which on $\Gamma_\mathrm{nm}$ corresponds to the points $a$ and $b$, and $\widetilde{\psi}$ be the  pullback of $\psi$ onto $\Gamma_\mathrm{nm}$. We have
\begin{equation}\label{eq1}
\widetilde{\psi}(a) = \widetilde{\psi}(b).
\end{equation}

\vskip4mm

{\bf Theorem B\cite{MT1}.}
{\sl For a singular spectral curve, theorem A remains true if we replace the genus $g$ with the arithmetic genus $g_a$, assuming that the points $P_1,\dots,P_n$ and $Q_1,\dots,Q_n$ are nonsingular points. }

\vskip4mm

In defining the divisor of the poles of the form $\Omega$, we do not include those singular points at which the form is regular. Let us explain that in more details.

A meromorphic form on a singular curve $\Gamma$ is defined as a meromorphic form on its
normalization $\Gamma_\mathrm{nm}$. The form is called regular at the point $P\in\Gamma$ if
either the point $P$ is not singular and in its preimage on $\Gamma_\mathrm{nm}$ the form $\Omega$ has no singularity, or $P$ is singular and
for all meromorphic functions on $\Gamma_\mathrm{nm}$, which descend to functions on
$\Gamma$ and have no poles in the preimages of $P$, the following relation holds:
$$
\sum_{\pi^{-1}(P)} \mathrm{Res} (f \Omega) = 0.
$$
For example, if $P$ is a double point, to which on $\Gamma_\mathrm{nm}$ there correspond the points
$a$ and $b$, then
the regularity condition at this point is written as
$$
\textrm{Res}_a \Omega + \textrm{Res}_b \Omega = 0.
$$

If for $1$-form (differential) $\Omega$ all points are regular, then we say that $\Omega$
is regular on $\Gamma$. The dimension of the space of regular
$1$-the forms is equal to the arithmetic genus of $\Gamma$.

\subsection{Curvilinear coordinates and torsion-free \\ sheaves of rank one on the spectral curves}

Let us consider our construction.

Let us assume that $\Gamma$ splits into irreducible nonsingular
components $\Gamma_j$ that intersect at double points.

The Baker-Akhiezer section $\psi$ is represented as a set of functions $\psi_j$ on the components.
These functions can be Baker-Achiezer functions with essential singularities or meromorphic functions.
At double points, these functions satisfy the relations
$$
\psi_j (a) = \lambda \psi_k(b), \ \ \lambda \ne 1.
$$
The coefficient $\lambda$ depends on a pair of points $a$ and $b$: $\lambda = \lambda_{ab}$, but in the sequel, for the sake of brevity, we will skip that when possible. 

The constructed function $\psi$ (on $\Gamma_\mathrm{nm}$)
defines a global section of some torsion-free sheaf $L$ over 
$\Gamma$. $L$ is obtained from $\Gamma_1\times{\mathbb C},\dots,\Gamma_m\times{\mathbb C}$ by gluing fibers over the intersection points of the components. Namely, if the components $\Gamma_p$ and $\Gamma_q$ intersect at the points $a\in\Gamma_p$ and $b\in \Gamma_q$, then the fibers
$a\times {\mathbb C}$ and $b\times {\mathbb C}$ are identified by using a suitable automorphism ${\mathbb C}$ (the multiplication by the constant $\lambda$).

\begin{theorem}
Let $\Gamma$ have a holomorphic involution $\sigma : \Gamma\rightarrow\Gamma$ such that 
$\sigma(\Gamma_s)=\Gamma_{\sigma(s)}$.
The points $P_j$ and $Q_j$ are fixed by $\sigma$, $\sigma (k_j) = - k_j$, and at the intersection points of the components of the curve $\Gamma$, we have
\begin{equation}
\label{eq2}
  \psi_p (a) = \lambda \psi_q (b),\qquad   \psi_{\sigma(p)} (\sigma (a)) = \lambda_{\sigma} \psi_{\sigma(q)} (\sigma(b)),
\end{equation}
where $\lambda, \lambda_{\sigma}$ are some non-zero constants.

Assume also that on $\Gamma$ there exists a meromorphic 1-form $\Omega$, given by meromorphic 1-forms $\Omega_j$ on the components $\Gamma_j$, the divisor of zeros of $\Omega$ has the form
$$
(\Omega)_{0}= P_1 +\dots +P_n + \gamma +\sigma\gamma,
$$
the divisor of poles of $\Omega$ contains the divisor
$$
Q_1 + \dots + Q_n + R + \sigma R
$$
and in addition $\Omega$ admits simple poles at the intersection points of the components and at these points we have
\begin{equation}\label{eq3}
 \lambda \lambda_{\sigma} {\rm Res}_a \Omega_p + {\rm Res}_b \Omega_q = 0,\ \
 \lambda \lambda_{\sigma} {\rm Res}_{\sigma(a)} \Omega_{\sigma(p)} + {\rm Res}_{\sigma(q)} \Omega_{\sigma(q)} = 0.
\end{equation}
Then if
$$
{\rm Res}_{Q_1} \Omega =\dots = {\rm Res}_{Q_n} \Omega=1,
$$
then
$$
\partial_{u^i} x^1\partial_{u^j} x^1 +\dots+\partial_{u^i} x^n\partial_{u^j} x^n = 0, \quad i\neq j,
$$
where
$$
x^j(u) = \psi (u,Q_j).
$$
\end{theorem}

{\sc Proof.}
Let us consider the form
$$
\omega_{ij} = \partial_{u^i} \psi (u,P) \partial_{u^j} \psi (u, \sigma (P)) \Omega.
$$
The form $\omega_{ij}$ is defined  on each component $\Gamma_k$ by the form
$$
\omega^k_{ij}= \partial_{u^i} \psi (u,P) \partial_{u^j} \psi (u, \sigma (P)) \Omega_k.
$$
The form $\omega_{ij}$ has simple poles at the points $Q_j$, as well as simple poles at the intersection points of the components of the curve $\Gamma,$
and
$$
{\rm Res}_{Q_s} \omega_{ij} = \partial_{u^i} x_s \partial_{u^j} x_s {\rm Res}_{Q_s} \Omega.
$$
Let, as above, the components $\Gamma_p$ and $\Gamma_q$ intersect at the points $a\in\Gamma_p$ and $b\in\Gamma_q$. Then
there are equalities
$$
{\rm Res}_{a} \omega^p_{ij} = \partial_{u^i} \psi_p(a) \partial_{u^j} \psi_{\sigma(p)} ( \sigma(a))
{\rm Res}_{a} \Omega_p,
$$
$$
{\rm Res}_{b} \omega^q_{ij} = \partial_{u^i} \psi_q (b)  \partial_{u^j} \psi_{\sigma(q)} ( \sigma (b) )
{\rm Res}_{b} \Omega_q.
$$

Note that, by (\ref{eq2}) and (\ref{eq3}),
$$
 {\rm Res}_{a} \omega^p_{ij} + {\rm Res}_{b} \omega^q_{ij} = \lambda \partial_{u^i} \psi_q (b) \lambda_{\sigma} \partial_{u^j} \psi_{\sigma(q)} ( \sigma (b) ) {\rm Res}_{a} \Omega_p +
$$
$$
 \partial_{u^i} \psi_q (b) \partial_{u^j} \psi_{\sigma(q)} ( \sigma (b) ) {\rm Res}_{b} \Omega_q
$$
$$= \partial_{u^i} \psi_q (b) \partial_{u^j} \psi_{\sigma(q)} ( \sigma (b) ) \left( \lambda \lambda_{\sigma} {\rm Res}_{a} \Omega_p +{\rm Res}_{b} \Omega_q \right) = 0.
$$

Then the sum of all residues of $\omega_{ij}$ for all components and all poles is equal to $0$. This gives
$$
 {\rm Res}_{Q_1} \omega_{ij} + \dots + {\rm Res}_{Q_n} \omega_{ij} = {\rm Res}_{Q_1} \Omega \partial_{u^i} x^1 \partial_{u^j} x^1 +\dots+ {\rm Res}_{Q_n} \Omega \partial_{u^i} x^n \partial_{u^j} x^n
$$
$$
 ={\rm Res}_{Q_1} \Omega(\partial_{u^i} x^1 \partial_{u^j} x^1 +\dots+  \partial_{u^i} x^n \partial_{u^j} x^n)= 0. \quad \qed
$$
Theorem is proved.

We describe the spectral data corresponding to the real-valued coordinates.

\begin{theorem}
Let there be an antiholomorphic involution on the spectral curve
 $\tau : \Gamma \rightarrow \Gamma$, $\tau:\Gamma_j=\Gamma_{\tau(j)}$ such that
 $$
   \tau(Q_j)=Q_j, \quad \tau(P_j)=P_j, \quad \tau(k_i^{-1}) = \overline{k_i^{-1}}, \quad \tau^{*} (\Omega) = \overline{\Omega},
 $$
$$
 \tau(\gamma) = \gamma, \quad \tau(R) = R.
$$
Therewith we assume that
$$
 \psi(R_j)=d_j,\quad \psi (\tau (R_j))=d_{\tau (j)},\quad \bar{d}_j=d_{\tau (j)},
$$
and at the intersection points of the components we have
$$
\psi_p (a) = \lambda \psi_q (b),\quad   \psi_{\tau(p)} (\tau (a)) = \lambda_{\tau} \psi_{\tau(q)} (\tau(b)),\quad
 \bar{\lambda}=\lambda_{\tau}.
$$
Then
$$
\psi(u,P) = \overline{\psi(u,\tau(P))},
$$
in particular, $x^j(u) = \psi(u,Q_j)$ are real-valued functions.
\end{theorem}

Theorem follows from the fact that the Baker-Akhiezer functions $\psi(u,P)$ and $\overline{\psi(u,\tau(P))}$ correspond to the same spectral data, therefore, by the uniqueness, they coincide.

{\sc Example 1.} Let $\Gamma$ be a reducible spectral curve consisting of two irreducible components $\Gamma_1$ and $\Gamma_2$
isomorphic to ${\mathbb C}P^1$ with affine coordinates $z_1$ and $z_2$ and the intersections at the points $a,-a\in\Gamma_1$ and $b,-b\in\Gamma_2$.

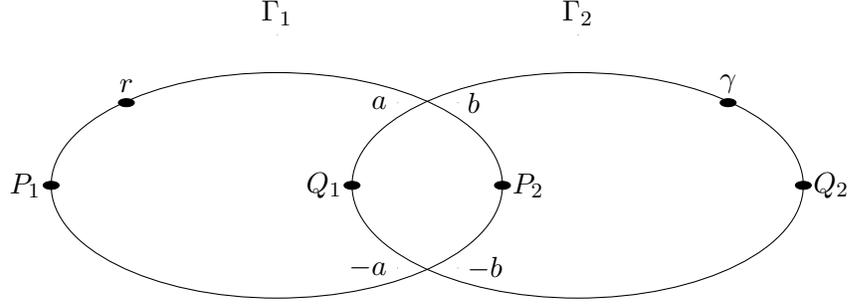
\begin{tikzpicture}[xscale=2]
     \draw (-1,0) circle (1.5);
     \draw (1,0) circle (1.5);
     \filldraw[black] (-2.5,0) coordinate (p1) circle (1.5pt) node[anchor=east] {$P_1$};
     \filldraw[black] (0.5,0) coordinate (p2) circle (1.5pt) node[anchor=west] {$P_2$};
     \filldraw[black] (-0.5,0) coordinate (q1) circle (1.5pt) node[anchor=east] {$Q_1$};
     \filldraw[black] (2.5,0) coordinate (q2) circle (1.5pt) node[anchor=west] {$Q_2$};
     \filldraw[black] (-2,1.1) coordinate (r) circle (1.5pt) node[anchor=south] {$r$};
     \filldraw[black] (2,1.1) coordinate (gamma) circle (1.5pt) node[anchor=south] {$\gamma$};
     \filldraw[black] (-1,2) coordinate (Gamma1) circle (0.0pt) node[anchor=south] {$\Gamma_1$};
     \filldraw[black] (1,2) coordinate (Gamma2) circle (0.0pt) node[anchor=south] {$\Gamma_2$};
     \filldraw[black] (-0.2,1.1) coordinate (a) circle (0.0pt) node[anchor=east] {$a$};
     \filldraw[black] (0.2,1.1) coordinate (b) circle (0.0pt) node[anchor=west] {$b$};
     \filldraw[black] (-0.2,-1.1) coordinate (mina) circle (0.0pt) node[anchor=east] {$-a$};
     \filldraw[black] (0.2,-1.1) coordinate (minb) circle (0.0pt) node[anchor=west] {$-b$};
\end{tikzpicture}

Let us take the following spectral data:
$$
 P_1=\infty,P_2=0, R=r\in\Gamma_1, \quad \ Q_1=\infty,Q_2=0, \gamma\in\Gamma_2.
$$
The curve $\Gamma$ admits the involution $\sigma:z_j=-z_j$.
We put
$$
\psi_{1} = e^{u^1z_{1} + u^2/z_{1} }f(u), \quad \psi_{2} = \left( g(u) + \dfrac{h(u)}{z_{2} - \gamma} \right).
$$
The conditions at the intersection points and the normalizing condition are:
$$
 \psi_{1}(a) =\lambda \psi_{2}(b), \quad  \psi_{1}(-a) =\mu \psi_{2}(-b), \quad \psi_{1}(r) = d.
$$
From that we obtain
$$
f(u)=d e^{-ru^1-u^2/r},\
$$
$$
 g(u)=\frac{de^{-\frac{(a+r)(aru^1+u^2)}{ar}}(\gamma(\lambda-\mu e^{2au^1+2u^2/a})+b(\lambda+\mu e^{2au^1+2u^2/a}))}{2b\lambda\mu},\
$$
$$
h(u)=-\frac{de^{-\frac{(a+r)(aru^1+u^2)}{ar}}(b^2-\gamma^2)(\lambda-\mu e^{2au^1+2u^2/a})}{2b\lambda\mu}.
$$
The differential $\Omega$ is given by the following differentials on $\Gamma_1$ and $\Gamma_2$:
$$
\Omega_1=\frac{sz_1dz_1}{(z_1^2-a^2)(z_1^2-r^2)},\quad\Omega_2=\frac{(z_2^2-\gamma^2)dz_2}{z_2(z_2^2-b^2)},
$$
with $s$ constant. The condition
$$
\textrm{Res}_{Q_1} \Omega_2 =\textrm{Res}_{Q_2} \Omega_2
$$
gives $b=i\gamma$. The conditions
$$
 \lambda\mu\textrm{Res}_{a} \Omega_1+ \textrm{Res}_{b} \Omega_2=0,
$$
$$
\lambda\mu\textrm{Res}_{-a} \Omega_1+ \textrm{Res}_{-b} \Omega_2=0,
$$
are satisfied for
$$
s=\frac{(r^2-a^2)(b^2-\gamma^2)}{b^2\lambda\mu}.
$$

We take 
$$
a=i,\quad d=1,\quad r=1,\quad\lambda=\bar{\mu}=\lambda_1+i\lambda_2.
$$
Then we derive the following curvilinear orthogonal coordinates on the plane:
$$
x=\psi_2(Q_1)=\frac{e^{-u^1-u^2}((\lambda_1-\lambda_2)\cos(u^1-u^2)+(\lambda_1+\lambda_2)\sin(u^1-u^2))}
{|\lambda|^2},
$$
$$
y=\psi_2(Q_2)=\frac{e^{-u^1-u^2}((\lambda_1+\lambda_2)\cos(u^1-u^2)-(\lambda_1-\lambda_2)\sin(u^1-u^2))}
{|\lambda|^2}.
$$

{\sc Example 2.}
Let the spectral curve look exactly the same as in Example 1. Take the following spectral data
$$
P=P_1+P_2,\quad Q=Q_1+Q_2,\quad R=R_1+R_2,\quad\gamma=\gamma_1+\gamma_2,
$$
$$
 P_1=\infty, \ P_2=0, \ R_1=r, \ R_2=-r, \ \gamma_1\in\Gamma_1,
 $$
 $$
 \ Q_1=\infty, \ Q_2=0, \ \gamma_2\in\Gamma_2.
$$
Let us put
$$
\psi_{1} = e^{u^1z_{1} + u^2/z_{1} }\left( f_1(u) + \frac{f_2(u)}{z_{1} - \gamma_1} \right),
\quad \psi_{2} = \left( g(u) + \frac{h(u)}{z_{2} - \gamma_2} \right).
$$

The conditions at the intersection points and the normalizing condition are
$$
\psi_{1}(a) =\lambda \psi_{2}(b), \quad  \psi_{1}(-a) =\mu \psi_{2}(-b), \quad \psi_{1}(r) = 1,\quad
\psi_{1}(-r) = 1.
$$
The differential $\Omega$ is given by the following differentials on $\Gamma_1$ and $\Gamma_2$:
$$
\Omega_1=\frac{sz_1(z_1^2-\gamma_1^2)dz_1}{(z_1^2-a^2)(z_1^2-r^2)^2},\quad\Omega_2=\frac{(z_2^2-\gamma_2^2)dz_2}{z_2(z_2^2-b^2)},
$$
with $s$ constant. The condition
$$
\textrm{Res}_{Q_1} \Omega_2 =\textrm{Res}_{Q_2} \Omega_2
$$
gives $b=i\gamma_2$. The conditions
$$
\lambda\mu\textrm{Res}_{a} \Omega_1+ \textrm{Res}_{b} \Omega_2=0,
$$
$$
\lambda\mu\textrm{Res}_{-a} \Omega_1+ \textrm{Res}_{-b} \Omega_2=0,
$$
are satisfied for
$$
s=-\frac{2(a^2-r^2)^2}{(a^2-\gamma_1^2)\lambda\mu}.
$$
We put 
$$
\lambda=\bar{\mu}=\lambda_1+i\lambda_2.
$$
To simplify the formulas, we also put
$$
a=i,r=i/2, \gamma_1=1.
$$
Then the coordinate functions take the form
$$
x=\psi_2(Q_1)=\frac{-(2\lambda_1+\lambda_2)\cos(\frac{3}{2}(u^1-2u^2))+
	3(2\lambda_1-\lambda_2)\cos(\frac{u^1}{2}+u^2)}{4|\lambda|^2}
$$
$$
+\frac{(\lambda_1-2\lambda_2)\sin(\frac{3}{2}(u^1-2u^2))+3(\lambda_1+2\lambda_2)
\sin(\frac{u^1}{2}+u^2)}{4|\lambda|^2} ,
$$

$$
y=\psi_2(Q_2)=\frac{(\lambda_1-2\lambda_2)\cos(\frac{3}{2}(u^1-2u^2))+
	3(\lambda_1+2\lambda_2)\cos(\frac{u^1}{2}+u^2)}{4|\lambda|^2}
$$
$$
+\frac{(2\lambda_1+\lambda_2)\sin(\frac{3}{2}(u^1-2u^2))+
	3(-2\lambda_1+\lambda_2)\sin(\frac{u^1}{2}+u^2)}{4|\lambda|^2}.
$$

\section{Darboux--Egorov metrics}

The diagonal metric
$$
ds^2=H_1^2(u)(du^1)^2+\dots+H_n^2(u)(du^n)^2
$$
in a domain of the Euclidean space is called the Darboux-Egorov metric if the rotation coefficients
$$
 \beta_{ij}=\frac{\partial_{u^i}H_j}{H_i},\ i\ne j
$$
are symmetric: $\beta_{ij}=\beta_{ji}$. Darboux-Egorov metrics play an important role in the theory of Frobenius manifolds \cite{D}. 

We indicate the restrictions on spectral data in Theorem 2, which lead to Darboux-Egorov metrics.
Let $\Gamma$ be a reducible spectral curve, each component of which be isomorphic to ${\mathbb C}P^1$ and $k_j=z_j$ be an affine coordinate on $\Gamma_j$, $j=1,\dots$.

Let the form $\Omega$ near $P_j$ have the following decomposition
$$
\Omega=\left(\frac{\varepsilon^2_j}{k_j}+O\left(\frac{1}{k^2_j}\right)\right)dk_j^{-1}.
$$

\begin{theorem}
Let each component $\Gamma_j, j=1,\dots,n$ contains points $P_j=\infty$ and $Q_j=0.$ Assume that each intersection point $q\in\Gamma_i\cap\Gamma_j$ has the same coordinates on each of the components $\Gamma_i$ and $\Gamma_j$:
\begin{equation}\label{eq4}
 z_i(q)=z_j(q).
\end{equation}
Then the flat metric has the form
\begin{equation}\label{eq5}
 ds^2=\varepsilon_1^2h_1^2(u)(du^1)^2+\dots+\varepsilon_n^2h_n^2(u)(du^n)^2
\end{equation}
and is the Darboux--Egorov metric.
\end{theorem}

{\sc Proof.} Let us show that in the coordinates $u^1,\dots,u^n$, the metric has the form (\ref{eq5}). To do this, we consider the form
$$
\omega_{ii}=\partial_{u^i}\psi(u,P)\partial_{u^i}\psi(u,\sigma(P))\Omega.
$$
This form has simple poles at the points $P_i,Q_1,\dots,Q_n$, as well as at the intersection points of the components of the spectral curve.
Let us sum up the residues of all forms defining $\omega_{ii}$ for all components. The residues corresponding to the intersection points of the components vanish (see proof of Theorem 1). We obtain
$$
(\partial_{u^j}x^1)^2+\dots+(\partial_{u^j}x^n)^2=\varepsilon_j^2h_j^2(u).
$$
Let us define on $\Gamma$ the function $f:\Gamma\rightarrow {\mathbb C}$ as follows. If $P\in\Gamma_j$, then we put $f(P)=z_j(P)$. By (\ref{eq4})  on the intersection points of the components, this function is correctly defined.

Next, consider the form
$$
\omega_{ij}=f(P)\frac{\partial_{u^i}\psi(u,z)}{h_i(u)}\frac{\partial_{u^j}\psi(u,\sigma(z))}{h_j(u)}\Omega.
$$
This form has simple poles at $P_i$ and $P_j$ and at the intersection points of the components. Summing up the sum of residues, we get:
$$
\frac{\partial_{u^j}h_i(u)\varepsilon_i^2}{h_j(u)}-\frac{\partial_{u^i}h_j(u)\varepsilon_j^2}{h_i(u)}=0,
$$
that is, $\beta_{ij}=\beta_{ji}$.
Theorem is proved.

{\sc Example 3.} Let $\Gamma$ be a reducible spectral curve, the same as in Example 1.
We assume that the components $\Gamma_1$ and $\Gamma_2$ intersect at points $a,-a\in\Gamma_1$, $a,-a\in\Gamma_2$.
We take the following spectral data:
$$
P_1=\infty, Q_1=0\in\Gamma_1, \quad \ P_2=\infty,Q_2=0,  R=r, \gamma\in\Gamma_2.
$$
The curve $\Gamma$ admits the involution $\sigma:z_j=-z_j.$
Let us put
$$
\psi_{1} = e^{u^1z_{1}}f(u), \quad \psi_{2} =  e^{u^2z_{2}}\left( g(u) + \dfrac{h(u)}{z_{2} - \gamma} \right).
$$
The conditions at the intersection points and the normalizing condition are
$$
\psi_{1}(a) =\lambda \psi_{2}(a), \quad  \psi_{1}(-a) =\mu \psi_{2}(-a), \quad \psi_{2}(r) = 1.
$$
The differential $\Omega$ is given by the following differentials on $\Gamma_1$ and $\Gamma_2$:
$$
\Omega_1=\frac{sdz_1}{z_1(z_1^2-a^2)},\quad\Omega_2=\frac{(z_2^2-\gamma^2)dz_2}{z_2(z_2^2-a^2)(z_2^2-r^2)},
$$
where $s$ is constant. The condition
$$
\textrm{Res}_{Q_1} \Omega_2 =\textrm{Res}_{Q_2} \Omega_2
$$
gives $s=\gamma^2/r^2$. The conditions
$$
\lambda\mu\textrm{Res}_{a} \Omega_1+ \textrm{Res}_{a} \Omega_2=0,
$$
$$
\lambda\mu\textrm{Res}_{-a} \Omega_1+ \textrm{Res}_{-a} \Omega_2=0,
$$
are satisfied for
$$
r=\frac{a\gamma\sqrt{\lambda\mu}}{\sqrt{\gamma^2(1+\lambda\mu)-a^2}}.
$$

Next, for simplicity of formulas, we put $a=1, \mu=5/(3\lambda), \gamma=2/3$. We obtain
$$
x^1=\frac{8e^{u^1-u^2}\lambda}{3(e^{2u^1}+e^{2u^2}\lambda^2)},\quad
  x^2=-\frac{2e^{2u^1-2u^2}-6\lambda^2}{3(e^{2u^1}+e^{2u^2}\lambda^2)}.
$$
Therewith the diagonal metric is a Darboux-Egorov metric
$$
ds^2=\frac{64e^{2u^1-2u^2}\lambda^2}{9(e^{2u^1}+e^{2u^2}\lambda^2)^2}(du^1)^2+
 \frac{16e^{-4u^2}(e^{2u^1}+3e^{2u^2}\lambda^2)^2}{9(e^{2u^1}+e^{2u^2}\lambda^2)^2}(du^2)^2.
$$

\section{Concluding remarks}

The spectral data for the polar and spherical coordinates were found in \cite{MT1} and thus
it was showed how they are embedded in the proposed scheme. {\it For elliptic coordinates on the plane, this problem remains open.}

This is related to the following problem.
By \cite[Lemma 3.1]{K1}, the function $\psi(u^1,\dots,u^n,Q)$ at each point $Q\in\Gamma$ of the spectral curve satisfies the equations
\begin{equation}
\label{second}
\frac{\partial^2 \psi}{\partial u^i \partial u^j} = \frac{\partial \log h_i}{\partial u^j}\frac{\partial \psi}{\partial u^i} +
\frac{\partial \log h_j}{\partial u^i}\frac{\partial \psi}{\partial u^j}, \ \ \ i \neq j,
\end{equation}
where
$$
ds^2 = \sum_i \varepsilon^2_i h_i^2(u) (du^i)^2, \ \ \ \varepsilon_i = \mathrm{const}
$$
and there is a flat metric for which $u^1,\dots,u^n$ are orthogonal curvilinear coordinates.
Since the coordinate $x^j(u)$ have the form $x^j(u) = \psi(u,Q_j), j=1,\dots,n$, then all coordinate functions must satisfy these equations.

For elliptic coordinates $(u^1,u^2)$ on the plane with the Euclidean coordinates $(x,y)$
we have
$$
x^2 = \frac{(u^1+a)(u^2+a)}{a-b}, \quad y^2 = \frac{(u^1+b)(u^2+b)}{b-a},
$$
where $a$ and $b$ are  constants such that $a>b$. The coordinates satisfy the inequalities
$$
-a < u^1 < -b < u^2,
$$
the metric takes the form
$$
dx^2 + dy^2 = \frac{1}{4}\frac{u^1-u^2}{(u^1+a)(u^1+b)} (du^1)^2 +
\frac{1}{4}\frac{u^2-u^1}{(u^2+a)(u^2+b)} (du^2)^2
$$
and the equation \eqref{second} reduces to the Euler--Poisson--Darboux equation $E(1/2,1/2)$:
$$
\frac{\partial^2 \psi}{\partial u^1 \partial u^2} = A\frac{1}{u^1-u^2}\frac{\partial \psi}{\partial u^1} -
B\frac{1}{u^1-u^2}\frac{\partial \psi}{\partial u^2}, \ \ A=B=\frac{1}{2}.
$$

The Euler-Poisson-Darboux equations $E(A,B)$ are well studied and a general solution has been found for them, depending on a pair of functional parameters \cite{Darboux}. To construct the elliptic coordinates by the finite-gap integration scheme, it is necessary to find in it a branch described by a Baker--Akhiezer function (or section).

In \cite{AVK} an integrable discretization of the Krichever scheme was proposed, consisting in constructing a similar Baker-Akhiezer function $\psi(k,l,Q)$, where integer parameters
$k,l\in {\mathbb Z}$ define a lattice on the plane. In this work, a discrete analogue of the Darboux-Egorov coordinates was distinguished.

In the mid-2010s, one of the authors (I.A.T.) proposed to construct discrete elliptic coordinates by a similar scheme, based on the equation $E(1/2,1/2)$. This was partially implemented
in \cite{BSST}, where some discretization of the equation $E(1/2,1/2)$ was used, unrelated to integrable systems. However, it would like {\it to construct discrete elliptic coordinates
in the framework of the theory of integrable systems}. Perhaps they will have very interesting properties, which is typical for exact solutions of integrable equations.

\end{document}